# Shear banding and cracking in unsaturated porous media through a nonlocal THM meshfree paradigm

Hossein Pashazad, Xiaoyu Song[✩]

*Engineering School of Sustainable Infrastructure and Environment*
*University of Florida, Gainesville, FL, USA*

**Abstract**

The mechanical behavior of unsaturated porous media under non-isothermal conditions plays a vital role in geo-hazards and geo-energy engineering (e.g., landslides triggered by fire and geothermal energy harvest and foundations). Temperature increase can trigger localized failure and cracking in unsaturated porous media. This article investigates the shear banding and cracking in unsaturated porous media under non-isothermal conditions through a thermo-hydro-mechanical (THM) periporomechanics (PPM) paradigm. PPM is a nonlocal formulation of classical poromechanics using integral equations, which is robust in simulating continuous and discontinuous deformation in porous media. As a new contribution, we formulate a nonlocal THM constitutive model for unsaturated porous media in the PPM paradigm in this study. The THM meshfree paradigm is implemented through an explicit Lagrangian meshfree algorithm. The return mapping algorithm is used to implement the nonlocal THM constitutive model numerically. Numerical examples are presented to assess the capability of the proposed THM mesh-free paradigm for modeling shear banding and cracking in unsaturated porous media under non-isothermal conditions. The numerical results are examined to study the effect of temperature variations on the formation of shear banding and cracking in unsaturated porous media.

## 1. Introduction

The thermo-hydro-mechanical (THM) behavior of unsaturated porous media, such as soils, plays a crucial role in various engineering applications, including nuclear waste disposal storage, pavement design, fault propagation, landslides, geothermal energy utilization, and the performance of buried high-voltage cables (e.g., [1–7]). Temperature variations can significantly impact the mechanical and physical properties of unsaturated soils, influencing parameters such as shear strength, deformation characteristics, fluid flow behavior, and mass transport properties at multiple length scales [8–16]. For instance, temperature changes can lead to complex behaviors in unsaturated soils, including volumetric strain or dilation, which may vary depending on factors like the overconsolidation ratio of the soil. Consequently, both physical experiments and numerical simulations are essential tools for investigating and understanding the coupled multi-physical processes involved in solid deformation, fluid flow, and heat conduction within thermally unsaturated soils (e.g., [9, 10, 17, 18]). These studies, such as the mesoscale finite element modeling of shear banding in thermal unsaturated soils [10] provide valuable insights into the behavior of unsaturated soils under the influence of temperature variations, contributing to more accurate and reliable engineering designs and assessments. In this study, as a new contribution, we develop a nonlocal mesh-free THM paradigm for modeling shear banding and cracking in unsaturated soils under elevated temperatures. For this purpose, we formulate a nonlocal THM constitutive model for unsaturated soils and implement the THM constitutive model into the meshfree periporomechanics (PPM) paradigm [19–28]. Next, we sequentially review the constitutive modeling of thermal unsaturated soils and the PPM paradigm.

---

[✩]Corresponding author
  *Email address:* xysong@ufl.edu (Xiaoyu Song)

Significant progress has been made in thermal constitutive modeling for unsaturated soils in recent decades, addressing the intricate interplay between thermal, mechanical, and hydraulic behavior under non-isothermal conditions [29–34]. These advancements are pivotal in understanding the response of unsaturated soils in a wide range of geotechnical applications. Numerous constitutive models have been formulated, each tailored to capture specific aspects of thermal-mechanical coupling in unsaturated soils (e.g., [17, 35, 36]). Some constitutive models have integrated thermal effects into established critical state theories [37], while others have explicitly accounted for temperature-induced alterations in the water retention curve [32]. Unified models have emerged, unifying both mechanical and thermal aspects, e.g., leveraging concepts from bounding surface theory [29, 36]. In addition, micro-structural-based constitutive models have been developed to elucidate the influence of temperature on capillary stress at solid-water-air interfaces [33]. Noteworthy contributions include hierarchical models that hierarchically incorporate hydro-mechanical hardening and thermal softening and models tailored to study cyclic behavior under varying thermal conditions [17, 36]. Collectively, these thermal constitutive models provide invaluable tools for comprehensively characterizing the behavior of unsaturated soils in response to temperature fluctuations, contributing to safer and more efficient engineering designs and geotechnical assessments. These advanced constitute models for thermal unsaturated soils have been implemented into the finite element program [38], which is robust for modeling continuous deformation in unsaturated soils but not for discontinuities such as shear bands and cracks. In the present study, we formulate a nonlocal thermal constitutive model for unsaturated soils and implement it into the mesh-free PPM paradigm to better study shear banding and cracking in thermal unsaturated soils.

PPM is a nonlocal formulation of classical coupled poromechanics, which is robust for modeling continuous and discontinuous mechanical and physical behavior of porous media [19–28, 39–42]. In PPM, equations of motion and mass balance are expressed as integral-differential equations [39, 40, 42]. PPM stands out for its natural ability to simulate multiphase discontinuities through field equations and material models [23, 26]. By using the stabilized multiphase correspondence principle, classical advanced constitutive models and physical laws are readily incorporated into PPM, enabling the modeling of coupled deformation, shear banding, and fracturing in porous media [19]. In PPM, the energy-based bond breakage criterion has been formulated for modeling cracks leveraging the effective force state concept [23, 26]. Furthermore, the large-deformation PPM through the updated Lagrangian framework was developed for unsaturated porous media in [27, 28]. The μPPM has been formulated to model dynamic shear bands and crack branching in porous media considering the rotational degree of freedom of the solid skeleton of porous media in [43–45]. In the present study, we investigate the shear banding and cracking in thermal unsaturated soils leveraging PPM. The PPM paradigm [19, 20] is used by incorporating a thermal constitutive model for unsaturated soils.

In this study, we delve into the intricate phenomena of shear banding and cracking within unsaturated porous media under non-isothermal conditions. A notable contribution of this study is the implementation of a classical THM material model tailored for unsaturated porous media into the computational meshfree PPM paradigm. A pivotal aspect of this integration is the utilization of a stabilized multiphase correspondence principle that effectively mitigates the zero-energy mode instability. Our implementation of the THM PPM paradigm is realized through an explicit Lagrangian meshfree algorithm. The return mapping algorithm in computational plasticity is used to numerically implement the THM constitutive model. To assess the capabilities of our proposed THM meshfree paradigm, we present numerical examples that illustrate its efficacy in modeling shear banding and cracking phenomena within unsaturated porous media under non-isothermal conditions. Our numerical results offer valuable insights into the intricate interplay between temperature variations and the formation of shear bands and cracks within unsaturated porous media.

The remainder of this article is organized as follows. Section 2 presents the mathematical formulation of the THM PPM framework including the thermal elastoplastic material model. Section 3 is dedicated to the numerical implementation of the proposed PPM paradigm. Section 4 presents numerical examples to assess the accuracy of the numerical implementation at the material point level and utilize the THM PPM paradigm to model dynamic shear banding and fracturing in unsaturated porous media under non-isothermal conditions, followed by a summary in Section 5. Throughout this work, we adopt the sign convention in continuum mechanics, wherein tensile forces and deformations under tension are considered positive. For pore fluid pressure, compression is positive, and tension is negative.

## 2. Mathematical formulation

In this section, we introduces the governing equation, the stabilized constitutive correspondence principle, the thermal elastoplastic material model, and the energy-based bond breakage criterion. In this study, we assume that the matric suction and temperature are known variables, i.e., one-way coupling. We also assume that no phase change exists between the three phases, i.e., solid, water and air.

### 2.1. Governing equation

In PPM, the porous media is represented by a set of mixed material points. A material point **X** has mechanical and physical interactions with any material point **X'** within its neighborhood, i.e., a spherical domain $H$ with a radius of $\delta$ called horizon. The bond $\underline{\xi}$ between material points **X** and **X'** is defined as $\underline{\xi} = $ **X' − X** in the reference configuration. For notation simplicity, the variables with no prime are associated with **X** and the variables with a prime means the variables associated with **X'**. For a partially saturated porous medium (i.e., solid, water, and air), assuming a weightless air phase, the total density is defined as

$$\rho = (1 - \phi)\rho_s + \phi S_r \rho_w, \tag{1}$$

where $\phi$ is the porosity, $\rho_s$ is the intrinsic density of the solid skeleton, $S_r$ is the degree of saturation, and $\rho_w$ is the intrinsic density of water. In this study, the degree of saturation is determined through a temperature-dependent water retention model for unsaturated soils at elevated temperatures [17], which reads

$$S_r = \left\{ \frac{1}{1 + \left[a_1 \gamma_\theta (\nu - 1)^{b_1} s\right]^n} \right\}^{-m}, \tag{2}$$

where $\nu$ is the specific volume of unsaturated porous media, $s$ is matric suction, $a_1$, $b_1$, $n$ and $m$ are material parameters, and $\gamma_\theta$ is a temperature-dependent air-entry matric suction. This variable can be determined by

$$\gamma_\theta = \left( \frac{a_2 + b_2 \theta_0}{a_2 + b_2 \theta} \right)^{b_2}, \tag{3}$$

where $a_2$ and $b_2$ are material parameters, $\theta_0$ is a reference temperature, and $\theta$ is the temperature of the mixture.

The motion equation for this porous medium in the PPM framework is written as

$$\rho \ddot{\boldsymbol{u}} = \int_{\mathcal{H}} \left( \underline{\mathcal{T}} - \underline{\mathcal{T}}' \right) d\mathcal{V}' - \rho \boldsymbol{g}, \tag{4}$$

where $\rho$ is the total density as defined in (1), $\ddot{u}$ is the acceleration, $\underline{T}$ and $\underline{T}'$ are the total force states (i.e., along the bond), and g is the gravity acceleration. Through the effective state concept [42] and assuming that matric suction and temperature are given, the motion equation for the porous media can degenerate into the motion equation for the solid phase as

$$\rho^s \ddot{\boldsymbol{u}} = \int_{\mathcal{H}} \left( \underline{\overline{\mathcal{T}}} - \underline{\overline{\mathcal{T}}}' \right) d\mathcal{V}' - \rho^s \boldsymbol{g}, \tag{5}$$

where $\rho_s$ is the partial density of the solid phase, and $\underline{T}$ and $\underline{T}'$ are the effective force states. Assuming the passive air pressure (i.e., zero air pressure), the effective force state at **X** is defined as

$$\overline{\underline{\mathscr{T}}} = \underline{\mathscr{T}} - S_r \overline{\underline{\mathscr{T}}}_w. \tag{6}$$

It is noted that the impact of temperature and matric suction on the mechanical behavior of unsaturated soils are considered through the thermal constitutive model given the temperature and matric suction. In what follows, we present the kinematics of the solid phase.

*2.2. Kinematics*

In PPM, the Lagrangian coordinate is used to modeling the solid phase [23]. Let y and **y'** be the positions of material points **X** and **X'** in the current configuration, respectively. Let u and **u'** be the displacements of material points **X** and **X'**, respectively. The deformation state and the displacement states are defined as

$$\underline{\mathscr{Y}} = \boldsymbol{y}' - \boldsymbol{y}, \tag{7}$$

$$\underline{\mathscr{U}} = \boldsymbol{u}' - \boldsymbol{u}. \tag{8}$$

Given $\underline{Y}$, the deformation gradient tensor in PPM [22] is defined as

$$\boldsymbol{F} = \left[ \int_{\mathscr{H}} \omega \left( \underline{\mathscr{Y}} \otimes \underline{\boldsymbol{\xi}} \right) d\mathscr{V}' \right] \mathscr{K}^{-1}, \tag{9}$$

where $\underline{\omega}$ is a weighting function, and $K$ is the shape tensor [42]. The shape tensor is defined as

$$\mathscr{K} = \int_{\mathscr{H}} \omega \left( \underline{\boldsymbol{\xi}} \otimes \underline{\boldsymbol{\xi}} \right) d\mathscr{V}'. \tag{10}$$

It is noted that the shape tensor $K$ is defined referred to the reference configuration. Then, it follows from (10), (9), and (8), the rate form of the deformation gradient tensor can be written as

$$\dot{\boldsymbol{F}} = \left[ \int_{\mathscr{H}} \omega \left( \underline{\dot{\mathscr{U}}} \otimes \underline{\boldsymbol{\xi}} \right) d\mathscr{V}' \right] \mathscr{K}^{-1}. \tag{11}$$

From (9) and (11), the velocity gradient tensor is determined as

$$\boldsymbol{L} = \dot{\boldsymbol{F}} \boldsymbol{F}^{-1}. \tag{12}$$

Given (12), the rate of deformation tensor D can be computed as

$$\boldsymbol{D} = \frac{1}{2} \left( \boldsymbol{L} + \boldsymbol{L}^T \right). \tag{13}$$

According to the polar decomposition theorem, the nonlocal deformation gradient F can be decomposed as

$$\boldsymbol{F} = \boldsymbol{R} \boldsymbol{U}, \tag{14}$$

where R is the rotation tensor that is a proper orthogonal tensor, and U is the right stretch tensor that is a symmetric positive-definite tensor. The unrotated rate of deformation tensor d can be obtained by

$$\hat{\boldsymbol{d}} = \boldsymbol{R} \boldsymbol{D} \boldsymbol{R}^T, \tag{15}$$

where the superscript T is the transpose operator.

Given the unrotated rate of deformation tensor, the strain increment can be written as

$$\Delta \varepsilon = \Delta t \hat{d}, \tag{16}$$

where $\Delta t$ is the time increment. Finally, given (9) the porosity [46] in the current configuration is written as

$$\phi = 1 - \frac{(1-\phi_0)}{J}, \tag{17}$$

where $J$ is the Jacobian of the nonlocal deformation gradient and $\phi_0$ is the initial porosity. We note that in this study the soil water retention curve is dependent on the porosity as introduced in Section 2.3.2. Next, we introduce the stabilized constitutive correspondence principle through which the advanced thermal constitutive model is implemented into the meshfree PPM paradigm.

*2.3. Correspondence THM constitutive model*

To complete (5), a constitutive model is needed to determine the effective force state. In this study, the stabilized constitutive correspondence principle [19] is used to implement an advanced thermal constitutive model for unsaturated soils.

*2.3.1. Constitutive correspondence principle*

The constitutive correspondence principle is based on the notion that the internal energy in a porous body from the local formulation in classical poromechanics is equal to that from the nonlocal formulation in periporomechanics. We refer to [19, 24, 42] for the detailed derivation. The effective force state in PPM can be written in terms of the effective Piola stress as

$$\overline{\underline{\mathscr{T}}} = \omega \overline{P} \mathscr{K}^{-1} \underline{\xi}, \tag{18}$$

where $\overline{P}$ is the effective Piola stress, which can be obtained from the local constitutive model given the nonlocal deformation gradient. It is note that, assuming passive air pressure (i.e., atmospheric air pressure), the effective stress $\overline{\sigma}$ is written as

$$\overline{\sigma} = \sigma - S_r p_w \mathbf{1}, \tag{19}$$

where $\sigma$ is the total Cauchy stress tensor, $p_w$ is pore water pressure, and $\mathbf{1}$ is the second-order identity tensor. Thus, it follows from (18), (6), and (19) that the fluid force state can be written as

$$\overline{\underline{\mathscr{T}}}_w = \omega p_w \mathbf{1} \mathscr{K}^{-1} \underline{\xi}. \tag{20}$$

Note that in (20) the small deformation of solid is assumed.

From (18) the effective force state can be computed from a thermal elasto-plastic constitutive model for unsaturated soils given matric suction, temperature change, and the nonlocal deformation gradient. The effective Piola stress can be written in terms of the unrotated Cauchy stress as

$$P = J\hat{\sigma} F^{-T}, \tag{21}$$

The unrotated effective Cauchy stress reads

$$\hat{\overline{\sigma}} = R\overline{\sigma} R^T, \tag{22}$$

where $\overline{\sigma}$ can be determined from an advanced thermal constitutive model for unsaturated soils. Next, we introduce the thermal elastoplastic model for unsaturated soils.

*2.3.2. Thermal elastoplastic model for unsaturated soils*

In this study, the thermal elastoplastic constitutive model is formulated based on the critical state soil mechanics. Following the small strain theory, the total strain is additively decomposed to elastic and plastic components as

$$\varepsilon = \varepsilon^e + \varepsilon^p, \qquad (23)$$

where $\varepsilon^e$ is the elastic strain tensor and $\varepsilon^p$ is the plastic strain tensor. For the thermal elastic model, the total elastic strain is assumed to consist of the mechanical elastic strain and the thermal elastic strain. Thus, the total elastic strain is additively decomposed into to mechanical and thermal parts as

$$\varepsilon^e = \varepsilon^e_{me} + \varepsilon^e_\theta, \qquad (24)$$

where $\varepsilon^e_{me}$ is the mechanical elastic strain and $\varepsilon^e_\theta$ is the thermal elastic strain. Given a temperature change, the thermal elastic strain is determined as

$$\varepsilon^e_\theta = \beta_\theta (\theta - \theta_0) \mathbf{1}, \qquad (25)$$

where $B_\theta$ is the volumetric thermal expansion coefficient, which is assumed as a constant in this study, $\theta$ is the temperature of soils, and $\theta_0$ is a reference temperature. Given the total elastic strain, the effective stress can be written through a linear thermal elastic model as

$$\overline{\sigma} = \mathbf{C} : \varepsilon^e, \qquad (26)$$

where C is the fourth-order linear elastic stiffness tensor that reads

$$C_{ijkl} = K\, \delta_{ij}\, \delta_{kl} + \mu\, (\delta_{ik}\delta_{jl} + \delta_{il}\delta_{jk} - \tfrac{2}{3}\, \delta_{ij}\, \delta_{kl}), \qquad (27)$$

where i, j, k, l = 1, 2, 3, K is the elastic bulk modulus, and μ is Poisson's ratio.

Next, we present the thermal plastic model. First, we define the effective mean stress $\overline{p}$ and the deviatoric stress q as

$$\overline{p} = \frac{1}{3}\mathrm{tr}(\overline{\sigma}), \qquad (28)$$

$$q = \sqrt{\frac{3}{2}} \|\overline{\sigma} - \overline{p}\mathbf{1}\|, \qquad (29)$$

Where || || is the norm of a tensor. Following the modified Cam-Clay model [17], the yield function is written as

$$f = \overline{p}^2 - \overline{p}p_c + \frac{q^2}{M^2}, \qquad (30)$$

where M is the slope of the critical state line and $p_c$ is the apparent preconsolidation pressure. In this study, the apparent preconsolidation pressure depends on the volumetric plastic strain, matric suction, and temperature changes [17]. Specifically, the apparent preconsolidation pressure reads

$$p_c = -\exp(\hat{a})\,(-p_{c,0})^{\hat{b}} \left[1 - \alpha_\theta \log\left(\frac{\theta}{\theta_0}\right)\right], \qquad (31)$$

where

$$\hat{a} = \frac{N(\hat{c} - 1)}{\widetilde{\lambda}\hat{c} - \widetilde{\kappa}}, \qquad (32)$$

$$\hat{b} = \frac{\widetilde{\lambda} - \widetilde{\kappa}}{\widetilde{\lambda}\hat{c} - \widetilde{\kappa}}, \qquad (33)$$

$$\hat{c} = 1 - c_1 [1 - \exp(c_2\zeta)], \qquad (34)$$

and $\zeta$ is a bonding variable related to water meniscus between grains, N is the specific volume of the soil under a unit saturated preconsolidaiton pressure, $c_1$ and $c_2$ are constants [47], $p_{c,o}$ is the apparent preconsolidation pressure at the reference temperature. It is noted that the parameter $c_\zeta$ is the ratio between the specific volume of the virgin compression curve in the partially saturated state to the corresponding specific volume in the fully saturated state. The bonding variable $\zeta$ [17] at the reference temperature (i.e., ambient temperature) is defined as

$$\zeta = (1 - S_r)\hat{f}(s), \tag{35}$$

where $(1-S_r)$ accounts for the number of water menisci per unit soil volume and $\hat{f}(s)$ is the stabilizing normal force exerted by a single water meniscus. The latter is written as

$$\hat{f}(s) = 1 + \frac{s/p_{\text{atm}}}{10.7 + 2.4\,(s/p_{atm})}, \tag{36}$$

where $p_{\text{atm}}$ is the atmospheric pressure.

Adopting the associative flow rule, the total THM plastic strain is written as

$$\dot{\varepsilon}^p = \dot{\lambda}\frac{\partial f}{\partial \boldsymbol{\sigma}}, \tag{37}$$

where $\dot{\lambda}$ is a plastic multiplier, which is determined by the consistency condition. Next, we introduce the energy-based bond breakage criterion.

### 2.4. Energy-based bond breakage criterion

In this study, the energy-based bond breakage criterion [24] is adopted to detect the bond breakage in the THM PPM framework. The effective force state is used to determine the deformation energy. Thus, the energy density in a bond is obtained as

$$\mathcal{W} = \int_0^t \left(\underline{\mathcal{T}} - \underline{\mathcal{T}}'\right) \dot{\underline{u}}\, dt, \tag{38}$$

where t is the load time. In PPM, the broken bond is modeled through the influence function at the constitutive model level. In this study, the influence function for bond $\leftarrow$ is defined as

$$\underline{\omega} = \begin{cases} 1 & \text{for } \mathcal{W} < \mathcal{W}_{cr}, \\ 0 & \text{for } \mathcal{W} \geq \mathcal{W}_{cr}, \end{cases} \tag{39}$$

where $\mathcal{W}_{cr}$ is the critical bond energy density. Following linear elastic fracture mechanics, the critical bond energy density can be calculated from the critical energy release rate as

$$\mathcal{W}_{cr} = \frac{4\mathcal{G}_{cr}}{\pi\delta^4}. \tag{40}$$

where $\mathcal{G}_{cr}$ is the critical energy per unit fracture area. In PPM, when a bond breaks, it will not sustain any mechanical load. The local damage parameter D at a material point is defined as

$$\mathcal{D} = 1 - \frac{\int_{\mathcal{H}} \underline{\omega}\, d\mathcal{V}'}{\int_{\mathcal{H}} d\mathcal{V}'}. \tag{41}$$

In this study, it is assumed that the crack initiates when $D > 0.5$ at a material point. In the following section, we present the numerical implementation of the THM PPM paradigm.

## 3. Numerical implementation

The THM PPM paradigm is implemented numerically through an explicit Newmark scheme [38, 48] in time and a Lagrangian meshfree method in space. The return mapping algorithm in computational plasticity is adopted for implementing the nonlocal thermal elasto-plastic constitutive model at the material point level. Algorithm 1 summarizes the global explicit meshfree numerical scheme and the local return mapping algorithm at the material point.

*3.1. Global integration in time*

In this part, we present the time integration of the governing equations at each material point. In this study, the explicit Newmark scheme is adopted. Let $u_n$, $\dot{u}_n$, and $\ddot{u}_n$ be the displacement, velocity, and acceleration vectors at time step n. The predictors of displacement and velocity in a general Newmark scheme read

$$\tilde{\dot{u}}_{n+1} = \dot{u}_n + (1-\beta_1)\Delta t \ddot{u}_n, \tag{42}$$

$$\tilde{u}_{n+1} = u_n + \Delta t \dot{u}_n + (1-2\beta_2)\Delta t^2 \ddot{u}_n, \tag{43}$$

where $\beta_1$ and $\beta_2$ are the numerical integration parameters. Given (42) and (43), the effective force state can be determined from the thermal elastoplastic constitutive model introduced in Section 2.3.2. Then, the acceleration at time step n + 1 is determined by

$$\ddot{u}_{n+1} = \mathcal{M}_{n+1}^{-1}\left(\widetilde{\mathcal{T}}_{n+1} - \mathcal{M}_{n+1}g\right), \tag{44}$$

where $\mathcal{M}_{n+1}$ is the mass of the solid at time step n + 1 and $\mathcal{T}_{n+1}$ is the effective force at time step n + 1. The two terms for a material point i are written as

$$\mathcal{M}_{n+1} = \rho_s (1-\phi_{n+1 i})\mathcal{V}_i, \tag{45}$$

$$\widetilde{\mathcal{T}}_{n+1} = \sum_{j=1}^{\mathcal{N}_i}\left(\widetilde{\underline{\mathcal{T}}}_{n+1,ij} - \widetilde{\underline{\mathcal{T}}}'_{n+1,ji}\right)\mathcal{V}_j\mathcal{V}_i, \tag{46}$$

where $\mathcal{N}_i$ is the number of neighbor material points of material point i. From (44), the displacement and velocity at time step n + 1 can be obtained as

$$\dot{u}_{n+1} = \tilde{\dot{u}}_{n+1} + \beta_1 \Delta t \ddot{u}_{n+1}, \tag{47}$$

$$u_{n+1} = \tilde{u}_{n+1} + \beta_2 \Delta t^2 \ddot{u}_{n+1}. \tag{48}$$

In this study, the explicit central difference solution scheme is adopted, i.e., $\beta_1 = 1/2$ and $\beta_2 = 0$. The energy balance check is used to ensure numerical stability of the algorithm in time. The internal energy, external energy, and kinetic energy of the system at time step n + 1 are written as

$$\mathscr{W}_{\text{int},n+1} = \mathscr{W}_{\text{int},n} + \frac{\Delta t}{2}\left(\dot{u}_n + \frac{\Delta t}{2}\ddot{u}_n\right)(\mathcal{T}_n + \mathcal{T}_{n+1}), \tag{49}$$

$$\mathscr{W}_{\text{ext},n+1} = \mathscr{W}_{\text{ext},n} + \frac{\Delta t}{2}\left(\dot{u}_n + \frac{\Delta t}{2}\ddot{u}_n\right)(\mathcal{M}_n g + \mathcal{M}_{n+1} g), \tag{50}$$

$$\mathscr{W}_{\text{kin},n+1} = \frac{1}{2}\dot{u}_{n+1}\mathcal{M}_{n+1}\dot{u}_{n+1}. \tag{51}$$

The energy conservation criterion requires

$$|\mathscr{W}_{\text{kin},n+1} - \mathscr{W}_{\text{ext},n+1} + \mathscr{W}_{\text{int},n+1}| \leq \hat{\varepsilon}\max(\mathscr{W}_{\text{kin},n+1}, \mathscr{W}_{\text{int},n+1}, \mathscr{W}_{\text{ext},n+1}), \tag{52}$$

where $\hat{\varepsilon}$ is a small tolerance on the order of $10^{-2}$ [48].

*3.2. Implementation of the material model*

This part deals with the numerical implementation of the thermal elasto-plastic model at the material point level through the return mapping algorithm (e.g., [10, 18]). First, we present the procedure for determining the strain increment at a material point i. Given (43), the deformation state on bond ij at time step n + 1 is written as

$$\widetilde{\underline{\mathcal{Y}}}_{n+1,ij} = \underline{\mathcal{Y}}_{n,ij} + \Delta\widetilde{\underline{\mathcal{Y}}}_{n+1,ij}. \tag{53}$$

Then, the nonlocal deformation gradient at material point i at time step n + 1 is computed by

$$\widetilde{F}_{n+1,i} = \left[\sum_{j=1}^{\mathcal{N}_i} \omega \left(\widetilde{\mathcal{Y}}_{n+1,ij} \otimes \underline{\xi}_{ij}\right) \mathcal{V}_j\right] \mathcal{K}_i^{-1}. \tag{54}$$

The spatial velocity gradient at material point i at time step n + 1 is written as

$$\widetilde{L}_{n+1,i} = \left[\left(\sum_{j=1}^{\mathcal{N}_i} \omega \dot{\widetilde{\mathcal{U}}}_{n+1,ij} \otimes \underline{\xi}_{ij} \mathcal{V}_j\right) \mathcal{K}_i^{-1}\right] \widetilde{F}_{n+1,i}^{-1}. \tag{55}$$

The rate of deformation tensor at time step n + 1 is written as

$$\widetilde{D}_{n+1,i} = \frac{1}{2}\left(\widetilde{L}_{n+1,i} + \widetilde{L}_{n+1,i}^T\right), \tag{56}$$

Given (56), the unrotated rate of deformation tensor $d_{n+1,i}$ can be written as

$$\widetilde{d}_{n+1,i} = R_{n+1,i} \widetilde{D}_{n+1,i} R_{n+1,i}^T. \tag{57}$$

where $R_{n+1,i}$ is rigid body rotation at material point i at time step n + 1. Then, the incremental strain tensor at material point i at time step n + 1 is computed as

$$\Delta \varepsilon_i = \Delta t \widetilde{d}_{n+1,i}. \tag{58}$$

Second, we present the procedure for updating the effective stress, given the increments of mechanical strain, temperature, and/or matric suction, through the return mapping algorithm. For brevity in notation, the subscript i of the material point is omitted in the following presentation. Let $s_n$, $\theta_n$ and $\varepsilon^e_n$
be the suction, temperature, and elastic strain, respectively, at material point
i at time step n. Let $\Delta \varepsilon_{me}$, $\Delta \theta$, and $\Delta s$ be incremental mechanical strain tensor, temperature, and matric suction from time steps n to n + 1. Here, we assume no return mapping on the suction and temperature [17, 47]. In this case, the matric suction and temperature at time step n + 1 can be written as

$$s_{n+1} = s_n + \Delta s, \tag{59}$$
$$\theta_{n+1} = \theta_n + \Delta \theta. \tag{60}$$

Given (60), the incremental thermal elastic strain $\Delta \varepsilon^e_{n+1}$ is defined as

$$\Delta \varepsilon^e_\theta = \frac{1}{3} \beta_\theta \Delta \theta \mathbf{1}. \tag{61}$$

By freezing plastic deformation, the trial elastic strain is written as

$$\varepsilon^{e,\text{tr}}_{n+1} = \varepsilon^e_n + \Delta \varepsilon + \Delta \varepsilon^e_\theta, \tag{62}$$

Then, the trial specific volume, degree of saturation, and bonding variable at time step n + 1 can be updated as

$$\nu_{n+1} = \nu_n \exp\left(1 + tr(\varepsilon^{e,\text{tr}}_{n+1})\right), \tag{63}$$

$$S_{r,n+1} = \left[\frac{1}{1 + \left[\widetilde{a}\gamma_\theta (\nu_{n+1} - 1)^{\widetilde{b}} s_{n+1}\right]^n}\right]^{-m}, \tag{64}$$

$$\zeta_{n+1} = (1 - S_{r,n+1}) \hat{f}(s_{n+1}). \tag{65}$$

The trial preconsolidatation pressure at time step n + 1 can be obtained from equation(31). To conduct the return mapping algorithm in the elastic strain space we define the unknown vector as

$$\boldsymbol{x}_{n+1} = \{\varepsilon_{v,n+1}^e, \varepsilon_{d,n+1}^e, \Delta\lambda\}^T, \tag{66}$$

where $\varepsilon_{v,n+1}^e$ is the elastic volume strain, $\varepsilon_{d,n+1}^e$ is the elastic deviatoric strain, and $\Delta\lambda$ is the plastic multiplier at time step n + 1. The residual vector is defined as

$$\boldsymbol{r}_{n+1} = \{r_{1,n+1}, r_{2,n+1}, r_{3,n+1}\}^T, \tag{67}$$

The elements of the residual vector are defined as

$$r_{1,n+1} = \varepsilon_{v,n+1}^e - \varepsilon_{v,n+1}^{e,\text{tr}} + \Delta\lambda \left(\frac{\partial f}{\partial \overline{p}}\right)_{n+1} \tag{68}$$

$$r_{2,n+1} = \varepsilon_{d,n+1}^e - \varepsilon_{d,n+1}^{e,\text{tr}} + \Delta\lambda \left(\frac{\partial f}{\partial q}\right)_{n+1} \tag{69}$$

$$r_{3,n+1} = f_{n+1}^{tr} \tag{70}$$

where $\varepsilon_{v,n+1}^{e,tr}$ is the trial elastic volume strain and $\varepsilon_{d,n+1}^{e,tr}$ is the trial elastic deviatoric strain at time step n + 1. The unknown vector x can be solved following the Newton's method as follows.

$$\delta \boldsymbol{x}_{n+1}^k = -\left(\frac{\partial \boldsymbol{r}}{\partial \boldsymbol{x}}\bigg|_{n+1}^k\right)^{-1} \boldsymbol{r}_{n+1}^k \tag{71}$$

$$\boldsymbol{x}_{n+1}^{k+1} = \boldsymbol{x}_{n+1}^k + \delta \boldsymbol{x}_{n+1}^k \tag{72}$$

where k is the iteration number. The tangent matrix in (71) reads

$$\frac{\partial \boldsymbol{r}}{\partial \boldsymbol{x}}\bigg|_{n+1}^k = \begin{bmatrix} \frac{\partial r_1}{\partial \varepsilon_v^e} & \frac{\partial r_1}{\partial \varepsilon_d^e} & \frac{\partial r_1}{\partial \Delta\lambda} \\ \frac{\partial r_2}{\partial \varepsilon_v^e} & \frac{\partial r_2}{\partial \varepsilon_d^e} & \frac{\partial r_2}{\partial \Delta\lambda} \\ \frac{\partial r_3}{\partial \varepsilon_v^e} & \frac{\partial r_3}{\partial \varepsilon_d^e} & \frac{\partial r_3}{\partial \Delta\lambda} \end{bmatrix}_{n+1}^k. \tag{73}$$

After solving the elastic strain, the effective stress at time step n + 1 can be updated through (26). The unrotated effective stress is

$$\overline{\boldsymbol{\sigma}}'_{n+1} = R_{n+1}^T \overline{\boldsymbol{\sigma}}_{n+1} R_{n+1}. \tag{74}$$

$$\overline{\mathscr{T}}_{n+1} = \omega \overline{\boldsymbol{P}}_{n+1} \mathscr{K}^{-1} \boldsymbol{\xi}. \tag{75}$$

From (21), the effective Piola stress at time step n + 1 can be computed. Then, the effective force state at time step n + 1 can be written as

**Algorithm 1** Summary of the numerical integration of the THM PPM paradigm

Given: $u_n, \dot{u}_n, \ddot{u}_n, \theta_n, s_n, \Delta\theta, \Delta s, \Delta t$ and compute: $u_{n+1}, \dot{u}_{n+1}, \ddot{u}_{n+1}, \theta_{n+1}, s_{n+1}$

1: Update time $t_{n+1} = t_n + \Delta t$
2: **while** $t_{n+1} \leq t_f$ **do**
3:    **for** all points **do**
4:       Compute the velocity predictor $\tilde{\dot{u}}_{n+1}$ using (42)
5:       Apply boundary conditions
6:       Compute displacement predictor $\tilde{u}_{n+1}$ using (43)
7:       **for** each neighbor **do**
8:          Update deformation state $\underline{\mathscr{Y}}_{n+1}$ using (53)
9:          Compute deformation gradient tensor $F_{n+1}$ using (54)
10:      **end for**
11:      Compute unrotated rate of deformation tensor $d_{n+1}$ using (57)
12:      Update temperature $\theta_{n+1}$ using (60) and suction $s_{n+1}$ using (59)
13:      Update preconsolidation pressure $p_{c,n+1}$
14:      Compute trial elastic strain tensor $\varepsilon_{n+1}^{e,tr}$ using (62)
15:      Compute the trial effective stress $\bar{\sigma}_{n+1}^{tr}$
16:      Compute the trial yield function $f_{n+1}^{tr}$
17:      **if** $f_{n+1}^{tr} \leq 0$ **then**
18:         Update effective stress $\bar{\sigma}_{n+1} = \bar{\sigma}_{n+1}^{tr}$
19:      **else if** $f_{n+1}^{tr} > 0$ **then**
20:         Compute the residual $r_{n+1}^k$
21:         **if** $\|r_{n+1}^k\| \leq$ Tol **then**
22:            Go to line 30
23:         **else if** $\|r_{n+1}^k\| >$ Tol **then**
24:            Compute $(\partial r/\partial x)_{n+1}^k$ using (73)
25:            Solve $\delta x_{n+1}^k$ using (71)
26:            Update the $x_{n+1}^{k+1}$ using (72)
27:            $k \leftarrow k + 1$
28:            Go to line 20
29:         **end if**
30:         Update effective stress $\bar{\sigma}_{n+1}$ using (26)
31:      **end if**
32:      Compute the effective force state using (75)
33:      Compute $\mathcal{M}_{n+1}$ using (45)
34:      Solve acceleration $\ddot{u}_{n+1}$ using (44)
35:      Update velocity $\dot{u}_{n+1}$ using (47)
36:      Update displacement $u_{n+1}$ using (48)
37:      Compute kinematic energy $\mathscr{W}_{kin,n+1}$ using (51)
38:      Compute internal energy $\mathscr{W}_{int,n+1}$ using (49) and external energy $\mathscr{W}_{ext,n+1}$ using (50)
39:      Check energy balance
40:      **for** each neighbor **do**
41:         Compute bond energy $W$
42:         **if** $W > W_{cr}$ **then**
43:            Update influence function
44:            Update damage variable $\mathscr{D}_{n+1}$
45:         **end if**
46:      **end for**
47:    **end for**
48: **end while**
49: $n \leftarrow n + 1$

## 4. Numerical examples

In this section, we present three numerical examples to showcase the effectiveness of the THM PPM paradigm in modeling shear banding and cracking in unsaturated porous media under THM conditions. Example 1 focuses on the isoerror map to assess the accuracy of the proposed re-turn mapping algorithm at the material point level. Example 2 addresses shear banding in an unsaturated elasto-plastic porous material under biaxial compression and varying temperature conditions. Example 3 examines crack formation in a disk specimen of an unsaturated elastic porous material under displacement control loading with increasing temperature.

*4.1. Accuracy assessment with isoerror maps*

This example evaluates the precision of the return mapping algorithm at the material point level through numerical testing. To gauge the accuracy of our proposed implicit algorithm, we employ isoerror maps [49]. The relative error is defined as follows:

$$\text{Error} = \frac{\sqrt{(\boldsymbol{\sigma} - \boldsymbol{\sigma}^\star) : (\boldsymbol{\sigma} - \boldsymbol{\sigma}^\star)}}{\sqrt{\boldsymbol{\sigma}^\star : \boldsymbol{\sigma}^\star}} \times 100, \tag{76}$$

where Error represents the algorithm's output and $\sigma^*$ denotes the exact solution, determined for specific strain and temperature increments. Following the methodology in [49], the exact solution is attained by repeatedly subdividing increments until further division yields negligible changes in the numerical result. It is important to note, as pointed out in [49], that while this approach effectively evaluates the algorithm's overall accuracy, it is not a substitute for a comprehensive analysis of accuracy and stability [49]. For this numerical test, the input material parameters [17, 47] are as follows: bulk modulus K = 83 MPa, shear modulus μ = 18 MPa, reference pressure $p_{c,o}$ = -35 kPa, reference specific volume 1.9, elastic thermal coefficient 6.67 × $10^{-4}$, swelling/recompression index 0.03, compression index 0.11, critical state line slope M = 1, plastic thermal parameter 0.23, $a_1$ = 0.038 kPa$^{-1}$, $b_1$ = 3.49, $a_2$ = 335°C, $b_2$ = 1, n = 0.718, m = 0.632, N = 2.76, initial temperature 25°C, $c_1$ = 0.185, and $c_2$ = 1.42.

We consider three distinct cases, each with specific initial conditions. For all cases, the initial effective isotropic compression stress is set uniformly at -150 kPa, and the preconsolidation pressure is established at -250 kPa. For Case 1, the initial temperature is 25 °C, with a constant matric suction of 50 kPa. For Cases 2 and 3, the initial temperature is raised to 50 °C, and the constant matric suction is increased to 100 kPa. For all three cases, the maximum temperature increment is set at 10 °C, and the maximum volumetric strain increment is -2%. To visualize the accuracy, we utilize isoerror maps plotted on a plane defined by the volumetric strain increment and the temperature increment. These maps employ a color bar to represent the error percentage. Figure 1 displays the isoerror maps for Case 1. Figure 2 illustrates the isoerror maps for Case 2. Figure 3 shows the isoerror maps for Case 3. The results in Figures 1 - 3 indicate that greater algorithmic accuracy can be achieved by adopting smaller increments in both temperature and strain.

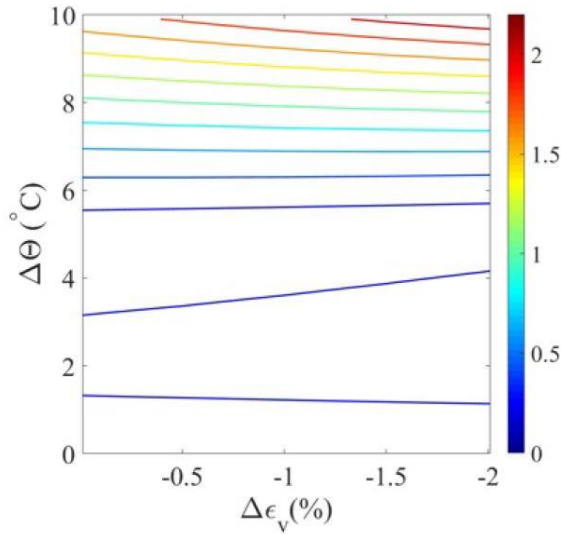

Figure 1: Isoerror map for case 1.

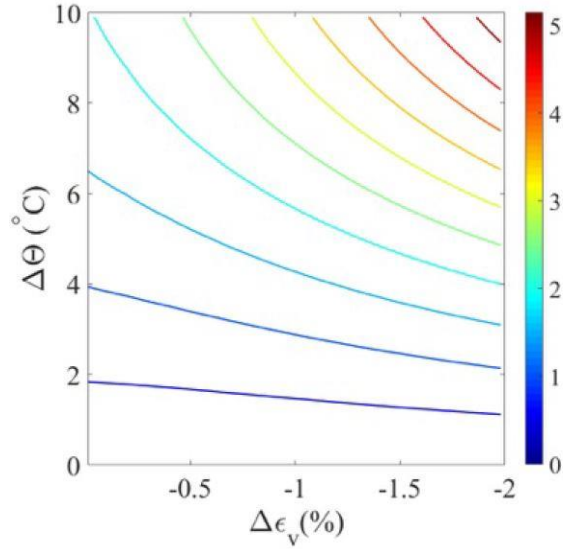

Figure 2: Isoerror map for case 2.

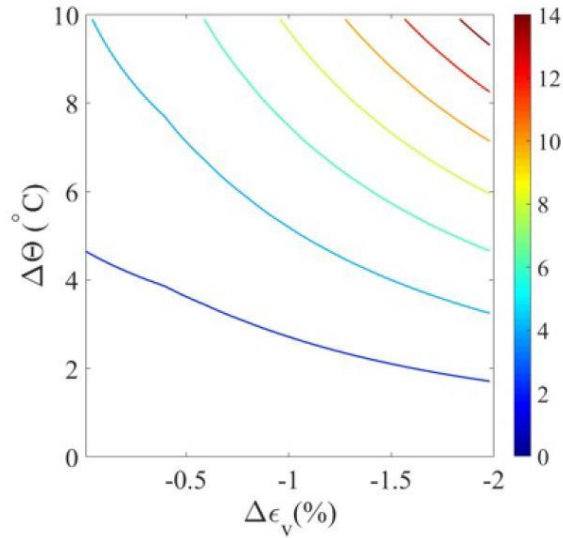

Figure 3: Isoerror map for case 3.

*4.2. Shear banding under non-isothermal conditions*

This example deals with the shear banding in thermal unsaturated porous media under dynamic loading conditions. Specifically, we investigate the influence of the effects of temperature and matric suction on shear banding. Figure 4 illustrates the model setup for this example. A vertical displacement of $u_y = 10$ mm is applied to the top boundary at a rate of 50 mm/s. A constant lateral confining pressure of 35 kPa is enforced on the left and right boundaries. The thermal elastoplastic constitutive model is utilized for this example. The input material parameters for the base simulation are: solid phase density 2000 kg/m³, bulk modulus $K = 83$ MPa, shear modulus $\mu = 18$ MPa, , elastic thermal coefficient $6.67 \times 10^{-4}$, reference pressure $p_{c,0} = -20$ kPa, reference specific volume 1.9, swelling index 0.03, compression index 0.11, critical state line slope $M = 1$, plastic thermal parameter $-0.23$, $a_1 = 0.038$ kPa, $b_1 = 3.49$, $a_2 = -335$ °C, $b_2 = 1$, $n = 0.718$, $m^o = 0.632$, $N = 2.76$, initial temperature 25 °C, $c_1 = 0.185$, and $c_2 = 1.42$. The specimen is discretized into a grid of $25 \times 50$ material points using a uniform grid spacing of 4 mm. The horizon is set to 8 mm, and the time increment is $1 \times 10^{-4}$ s.

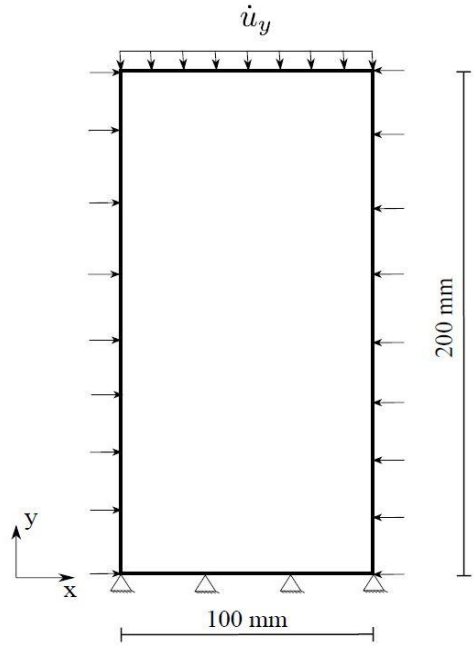

Figure 4: Model setup for the example of shear banding.

For the base simulation, a constant temperature of 25°C is prescribed within the problem domain. The matric suction decreases from 25 kPa to 10 kPa at a rate of 7.5 kPa/s. The results of the base simulation are presented in Figures 5 - 8. Figure 5 plots the loading curve on the top boundary, demonstrating a softening stage after the peak load due to the reduction of matric suction. Figure 6 displays the curve of deviatoric stress with vertical strain and the stress path (in the p  q space) of the point at the specimen center. The results indicate that the deviator stress increases with mean stress until it reaches the critical state line, after which it starts to decrease due to softening. Figure 7 presents snapshots of the equivalent plastic shear strain in the deformed configuration at three loading stages. Figure 8 provides snapshots of the plastic volumetric strain at the same three loading stages. It is important to note that a magnification factor of 5 is applied to all contours in this example. The results in Figures 7 and 8 demonstrate the development of two conjugate shear bands originating from the specimen center. Notably, in our nonlocal PPM framework, the initiation of shear banding does not require a weak element, as typically seen in finite element modeling of shear banding. Figure 8 shows that the plastic volumetric strain in the shear zone is positive, indicating dilatation.

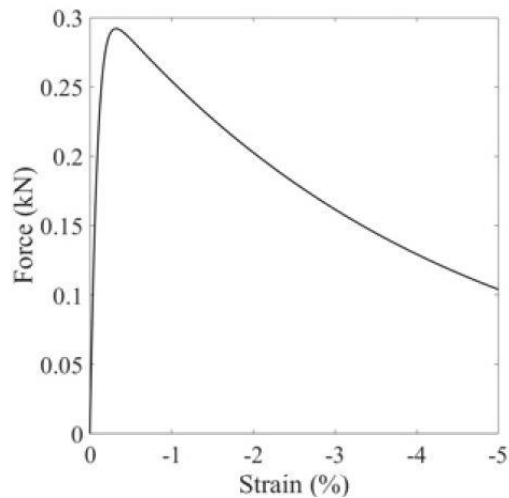

Figure 5: Loading curve on top boundary.

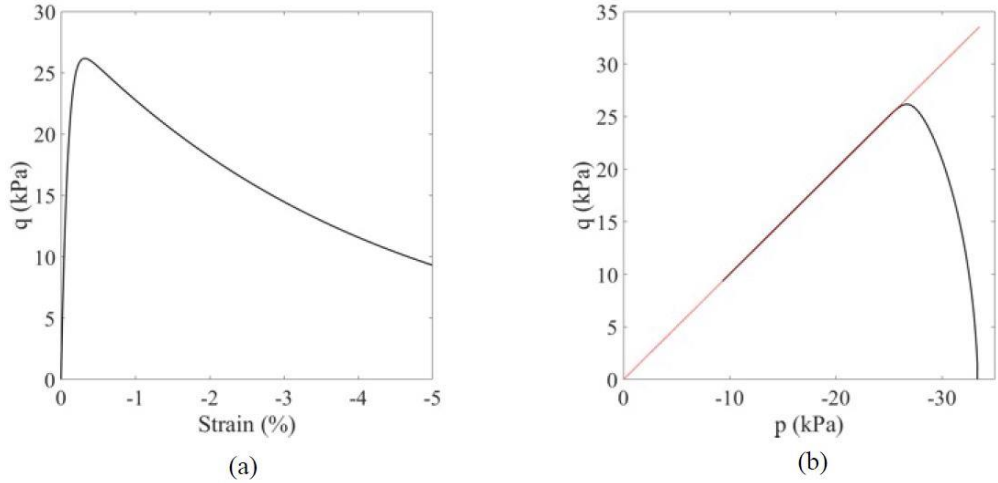

Figure 6: (a) Plot of the deviator stress versus vertical strain and (b) stress path in the p - q space at the point of the the specimen center. Note: The dash line in red is the critical state line).

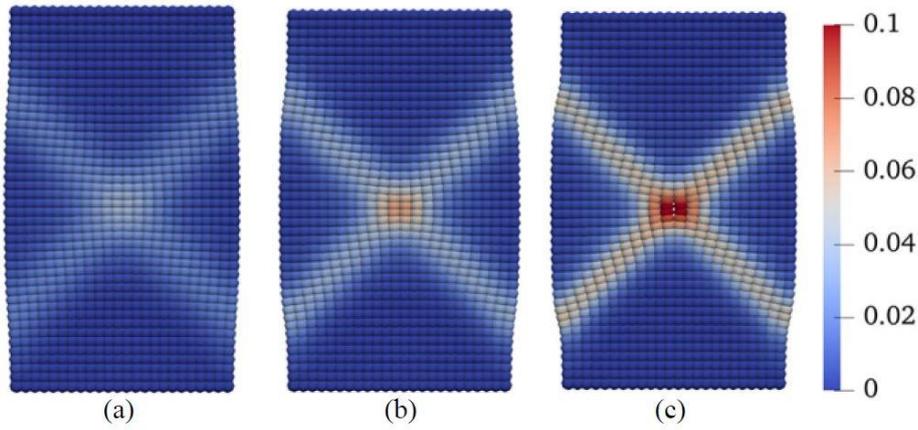

Figure 7: Contours of the equivalent plastic shear strain on the deformed configuration at (a) $u_y = 4$ mm, (b) $u_y = 7$ mm, and (c) $u_y = 10$ mm

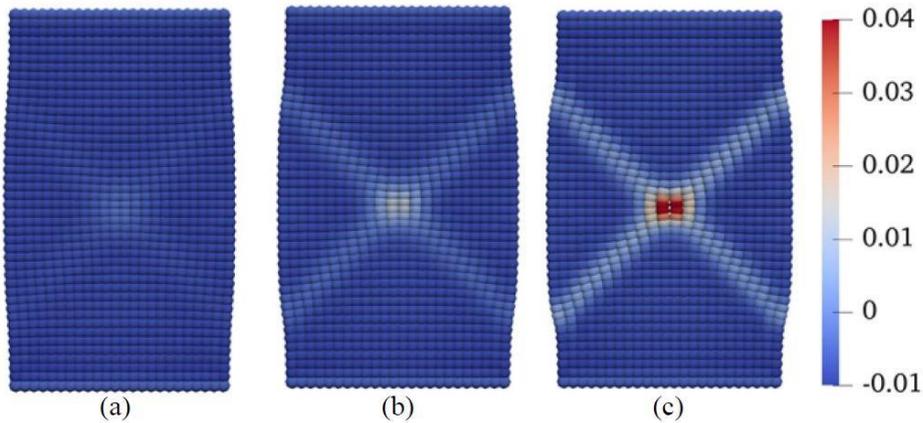

Figure 8: Contours of the plastic volume strain on the deformed configuration at (a) $u_y = 4$ mm, (b) $u_y = 7$ mm, and (c) $u_y = 10$ mm.

To investigate the impact of spatial discretization on the results, we examine two different spatial discretization schemes: one with a grid of $25 \times 50$ points and $\Delta x = 4$ mm (referred to as grid 1), and the other with a grid of $40 \times 80$ points and $\Delta x = 2.5$ mm (referred to as grid 2). Both simulations utilize the same horizon value of 8 mm, while all other conditions and parameters remain consistent with the base simulation. Figure 9 presents a comparison of the loading curves obtained from the two simulations. These two loading curves are identical until the onset of the softening stage. Figure 10 displays the contours of equivalent plastic shear strain at $u_y = 10$ mm for both simulations, while Figure 11 shows the contours of plastic volumetric strains at the same displacement level. The results from Figures 10 and 8 suggest that the choice of spatial discretization has a relatively minor influence on shear band formation, primarily due to the adoption of the same nonlocal length scale. In the subsequent sections, we investigate the impact of temperature on shear banding in unsaturated porous media at elevated temperatures. Three scenarios are considered: (i) Elevated constant temperature (scenario 1), (ii) Increasing temperature at constant suction (scenario 2), and (iii) Increasing temperature under decreasing suction (scenario 3).

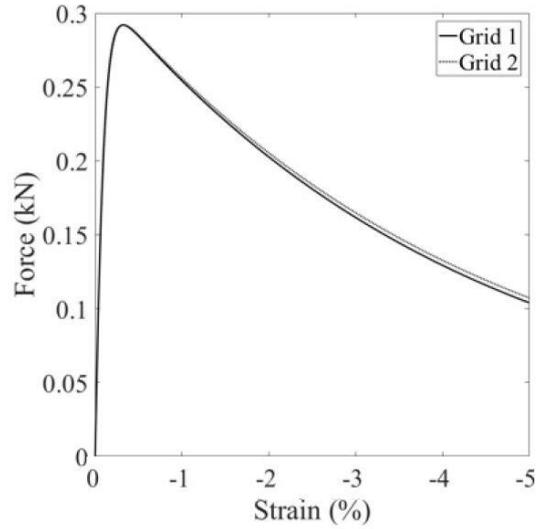

Figure 9: Comparison of the loading curves on the top boundary from the simulations with two spatial discretizations.

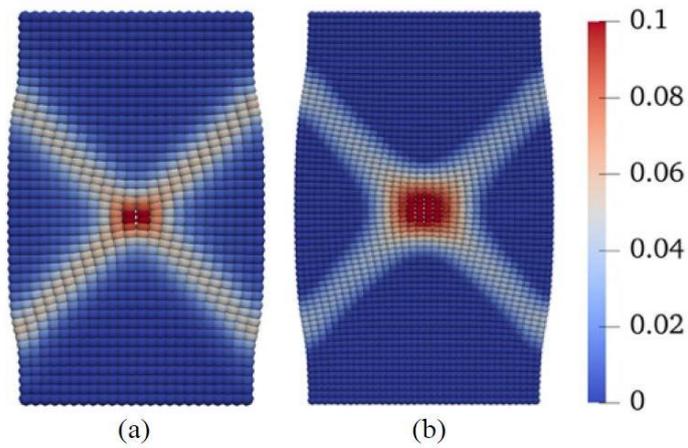

Figure 10: Contours of the equivalent plastic shear strain on the deformed configuration at $u_y = 10$ mm from the simulations with (a) grid 1 and (b) grid 2.

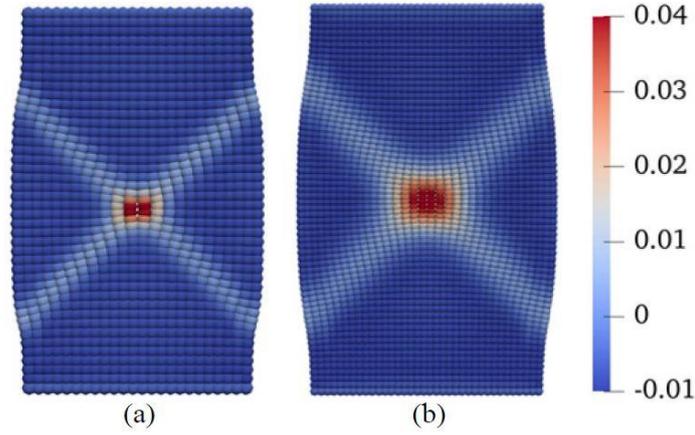

Figure 11: Contours of the plastic volume strain on the deformed configuration at $u_y = 10$ mm from the simulations with (a) grid 1 and (b) grid 2.

*4.2.1. Scenario 1: Elevated constant temperature*

In this scenario, we explore the influence of temperature on shear banding under constant suction conditions. To achieve this, we conduct numerical simulations at three different temperatures: 25°C, 50°C, and 75°C, all while maintaining a constant matric suction of 25 kPa. All other conditions and parameters remain consistent with the base simulation. The results of these simulations are presented in Figures 12 -15. Figure 12 compares the loading curves obtained from the three simulations. As shown in Figure 12, the loading capacity of the specimen decreases at higher temperatures due to the temperature-induced softening effect. Figure 13 illustrates the curves of deviator stress versus vertical strain and the stress paths at the specimen center from the three simulations. Regardless of the temperature, the stress paths at the same point demonstrate that the soil element reaches the same critical state line under loading. It's noteworthy that the critical state line remains consistent due to the adoption of the Bishop-type effective stress model for unsaturated soils. Figure 14 displays the contours of equivalent plastic shear strain at $u_y = 10$ mm from the three simulations, while Figure 15 presents the contours of plastic volumetric strain at the same displacement level. These results imply that temperature affects the magnitudes of dilation and shear strain within the specimen under the same mechanical load. Specifically, higher temperatures lead to more significant dilation and shear strain compared to lower temperatures. In the subsequent section, we further investigate the impact of varying temperature on the formation of shear banding.

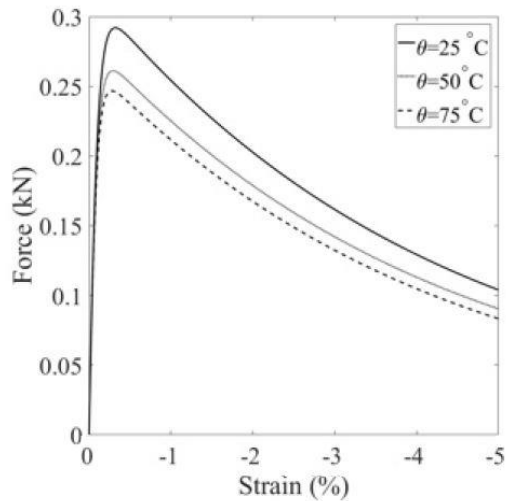

Figure 12: Comparing of the loading curves from the simulations under the three temperatures.

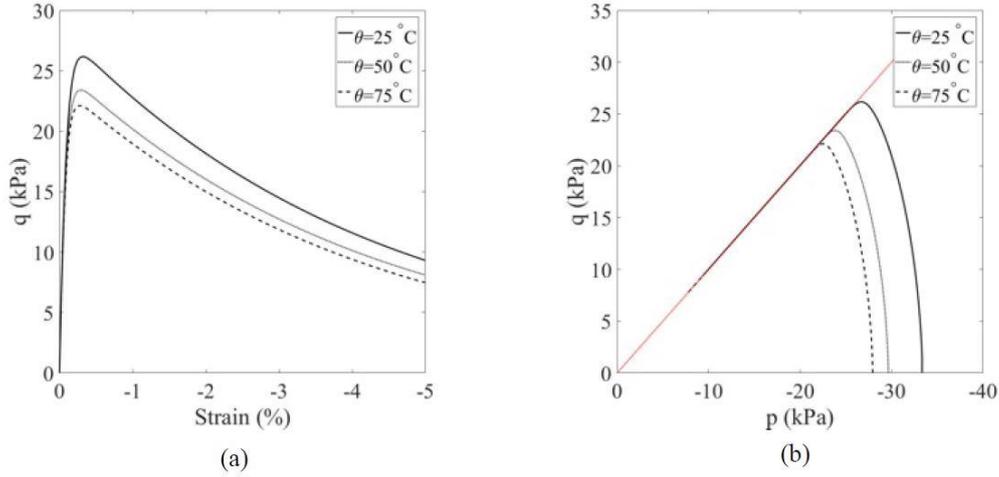

Figure 13: (a) Plot of deviatoric stress versus vertical strain, and (b) the stress paths in the p - q space at the specimen center.

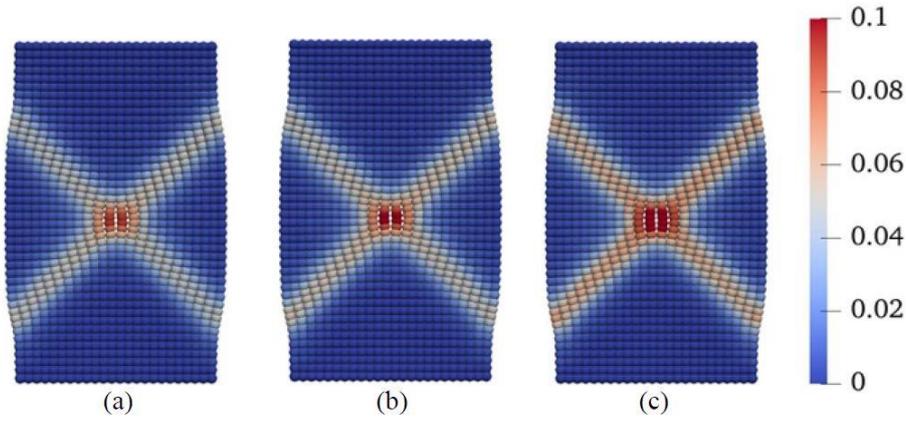

Figure 14: Contours of equivalent plastic shear strain on deformed configuration: (a) $\theta = 25°C$, (b) $\theta = 50°C$, and (c) $\theta = 75°C$.

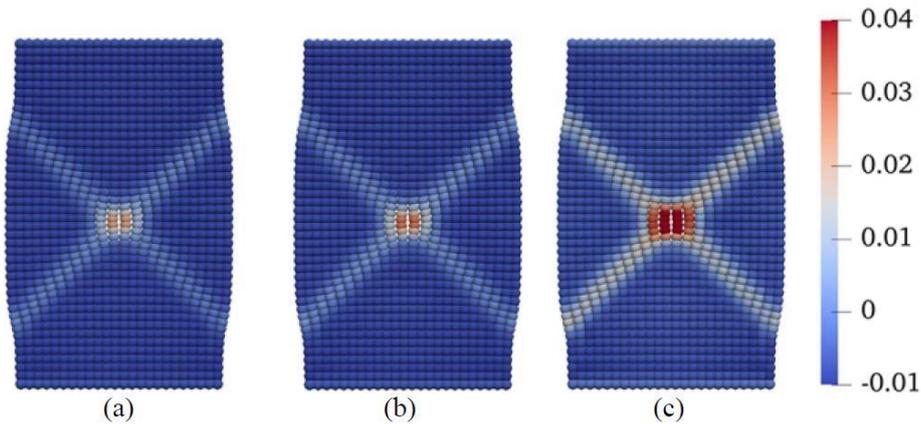

Figure 15: Contours of plastic volume strain on deformed configuration at $u_y = 10$ mm: (a) $\theta = 25°C$, (b) $\theta = 50°C$, and (c) $\theta = 75°C$.

*4.2.2. Scenario 2: Increasing temperature*

In this scenario, we examine the effect of temperature increase on the development of shear banding in unsaturated soils after the peak load. Specifically, we consider three different temperature changes applied after reaching the peak load of the base simulation: 25ºC, 25ºC, and 50ºC, all while maintaining a constant suction level of 25 kPa. It is important to note that the simulation with tl√ = 25ºC is included for comparison purposes. All other conditions and input parameters remain consistent with the base simulation. The results are presented in Figures 16 -

19. Figure 16 provides a comparison of the loading curves obtained from the three simulations. As depicted in Figure 16, increasing the temperature after the peak load has a notable effect on the post-localization regime in unsaturated soils, with a larger temperature increase resulting in a more significant reduction in strength. Figure 17 displays the curves of deviatoric stress versus vertical strain and the stress paths at the specimen center from the three simulations. The results in Figures 16 and 17 reinforce the influence of temperature increase on the post-localization behavior in unsaturated soils, where a larger temperature increase leads to a more pronounced strength reduction. Figure 18 presents the contours of equivalent plastic shear strain at $u_y = 10$ mm from the three simulations, while Figure 19 plots the contours of plastic volumetric strain at the same displacement level. These results, as shown in Figures 18 and 19, illustrate that a larger temperature increase, under the same conditions, results in more significant dilation and shear strain within the shear banding zone. In the subsequent section, we delve into the combined effect of varying temperature and suction on shear banding.

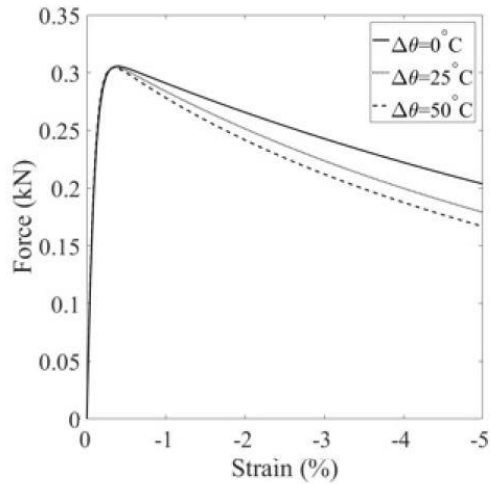

Figure 16: Comparison of the loading curves from the simulations with three temperature changes.

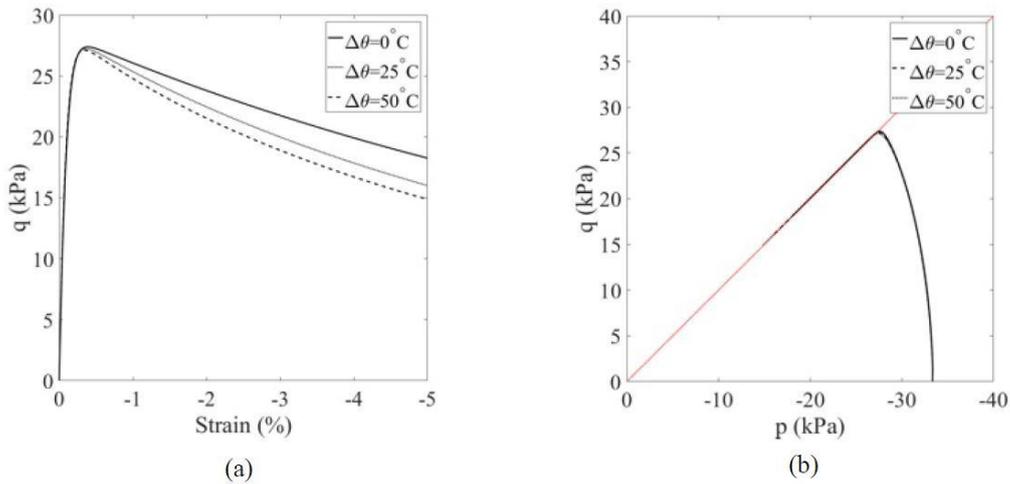

Figure 17: (a) Plot of deviatoric stress versus vertical strain, and (b) stress paths in the p-q space at the specimen center under three temperature changes.

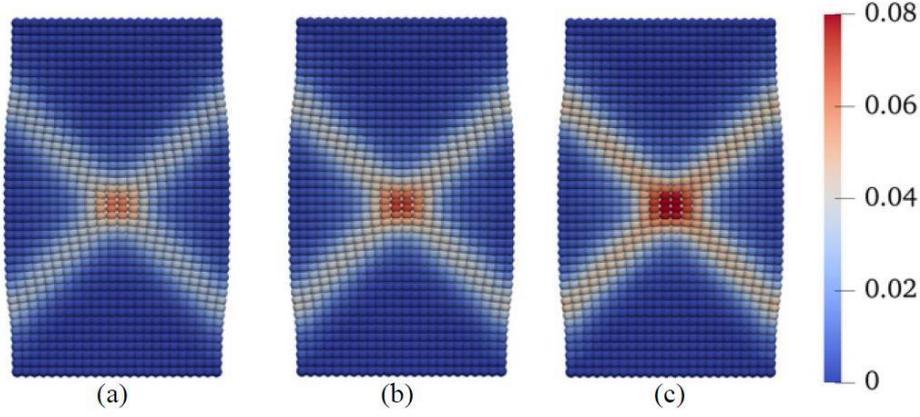

Figure 18: Contours of equivalent plastic shear strain on deformed configuration at $u_y = 10$ mm: (a) $\Delta\theta = 0°\text{C}$, (b) $\Delta\theta = 25°\text{C}$, and (c) $\Delta\theta = 50°\text{C}$.

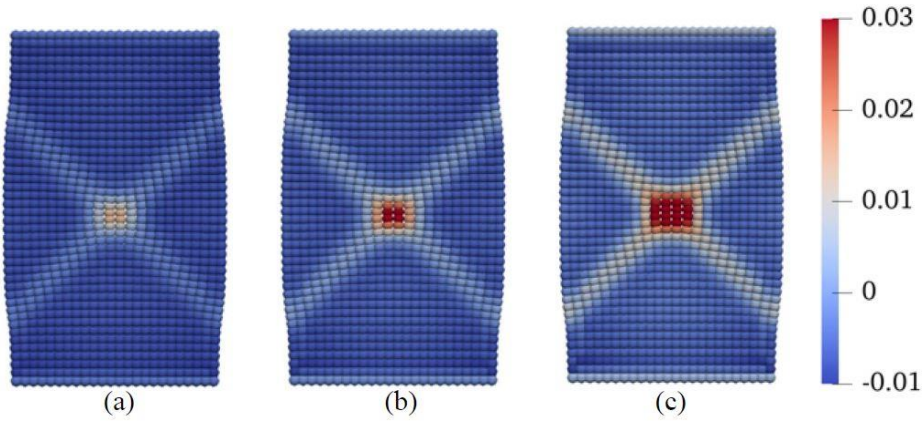

Figure 19: Contours of plastic volume strain on deformed configuration at $u_y = 10$ mm: (a) $\Delta\theta = 0°\text{C}$, (b) $\Delta\theta = 25°\text{C}$, and (c) $\Delta\theta = 50°\text{C}$.

*4.2.3. Scenario 3: Increasing temperature and decreasing suction*

In this scenario, we investigate the combined effect of increasing temperature and decreasing suction on shear banding instability in unsaturated soils. To achieve this, we consider three different temperature changes: 0ºC, 25ºC, and 50ºC, while concurrently decreasing suction from 25 kPa to 10 kPa. All other parameters and loading conditions are kept consistent with the base simulation, and the results are presented in Figures 20 - 23. Figure 20 displays the loading curves obtained from the three simulations. As shown in Figure 20, there is a notable reduction in loading capacity under the combined effect of temperature increase and suction reduction during the post-localization regime of unsaturated soils. Figure 21 presents the curves of deviatoric stress versus vertical strain and the stress paths from the three simulations. These results, depicted in Figure 21, further emphasize the impact of the coupling effect, showing that the soil reaches at critical state line at the same point for all three temperature-suction scenarios. Figure 22 compares the contours of equivalent plastic shear strain at $u_y$ = 10 mm in the deformed configuration for the three simulations. Meanwhile, Figure 23 compares the contours of plastic volumetric strain at $u_y$ = 10 mm. The results in Figures 20 - 23 illustrate a significant reduction in loading capacity under the coupling effect of temperature increase and suction reduction during the post-localization regime in unsaturated soils. In summary, these findings highlight the complex interplay between temperature and suction on shear banding instability in unsaturated soils.

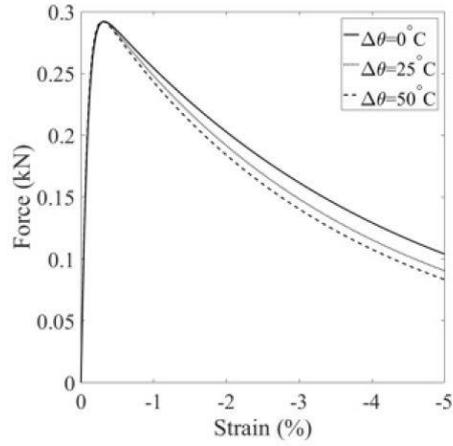

Figure 20: Comparison of the loading curves on the top boundary for Scenario 3.

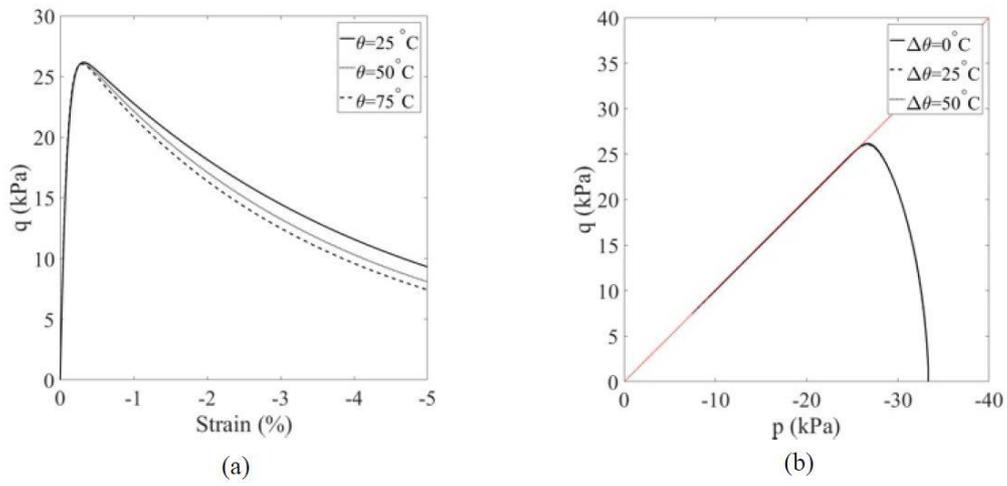

Figure 21: (a) Plot of deviatoric stress versus vertical strain, and (b) stress loading path in p versus q space (red dash line represents the critical state line) at the specimen center for Scenario 3.

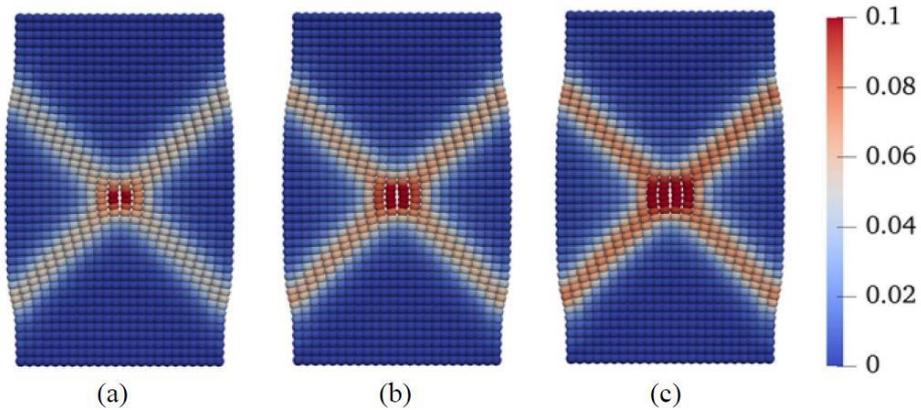

Figure 22: Contours of equivalent plastic shear strain on deformed configuration at $u_y = 10$ mm for Scenario 3: (a) $\Delta\theta = 0°$C, (b) $\Delta\theta = 25°$C, and (c) $\Delta\theta = 50°$C.

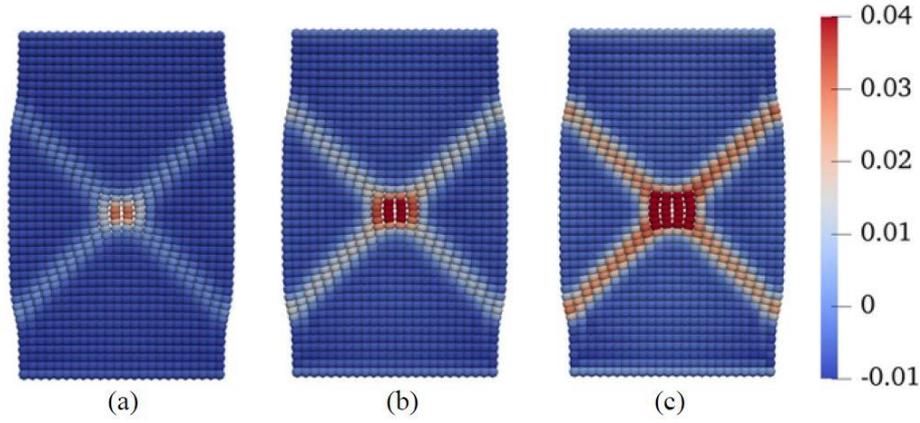

Figure 23: Contours of plastic volume strain on deformed configuration at $u_y = 10$ mm for Scenario 3: (a) $\Delta\theta = 0°C$, (b) $\Delta\theta = 25°C$, and (c) $\Delta\theta = 50°C$.

*4.3. Cracking in an elastic unsaturated disk specimen*

In this example, we focus on cracking phenomena in unsaturated elastic porous materials. The modeling of cracking is based on an energy-based bond breakage criterion. Specifically, we simulate cracking in a disk specimen. Figure 24 illustrates the disk specimen and its loading scheme. The disk has a radius of 200 mm and a thickness of 5 mm. As depicted in Figure 24, vertical displacement loads are applied to the top and bottom plates of the disk. The displacement load on each plate is set at u = 2.0 mm, with a loading rate of u˙ = 200 mm/s. The short-range forces within the PPM framework are employed to simulate the contact between the disk and the rigid plates [21]. The matric suction present in the specimen is s = 10 kPa. For this example, we adopt a thermo-elastic material model. The material parameters include: solid phase density of 2000 kg/m³, bulk modulus K = 83 MPa, shear modulus μ = 18 MPa.. In the base simulation scenario, the temperature is increased by 50ºC. The energy-based bond breakage criterion is implemented with a critical
energy release rate, $G_{cr}$ = 20 N/m. The specimen discretization involves 10408 points arranged in a uniform grid, with a spacing of 2.5 mm. The horizon size is set to δ = 4Δx. The simulation uses a time increment of $1 \times 10^{-5}$ s.

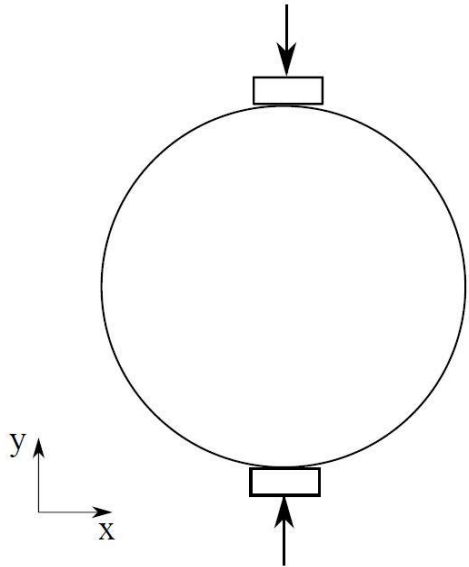

Figure 24: Model setup for the disk specimen under compression.

First, we present the results of the base simulation conducted under ambient temperature conditions. Figure 25 displays the loading curve applied to the top plate of the specimen, revealing a peak load of approximately 0.2 kN. To illustrate the progression of crack formation within the specimen, we provide contours of displacements at three distinct loading stages. Figure 26 depicts the vertical displacement contours on the deformed specimen at these stages, with a uniform magnification factor of 5 applied to all contours in this example. Similarly, Figure 27 presents the horizontal displacement contours under the same conditions. Analysis of the results shown in Figures 26 and 27 indicates the initiation of a crack at the specimen's center, which then extends towards the top and bottom of the disk. Notably, there is a discontinuity in the x-direction displacement along the vertical center line, whereas the vertical displacement remains continuous along the horizontal center line. This observation aligns with the expectations set by classical Brazilian testing, suggesting that the crack results from the discontinuous deformation in the x direction, as further evidenced by the deformed configurations depicted in Figures 26 and 27. Furthermore, Figure 28 illustrates the vertical normal stress contours on the deformed specimen at the three stages, highlighting the maximum compression stress occurring under the plate. Figure 29 focuses on the contours of horizontal normal stress elucidating that the specimen experiences tension in the x direction, with the crack process zone around the crack tip being under horizontal tension and vertical compression. Figure 30 shows the contours of the damage parameter at the same three stages. The results from Figures 28 to 30 collectively indicate that the crack formation is primarily due to tensile stress perpendicular to the specimen's vertical center line.

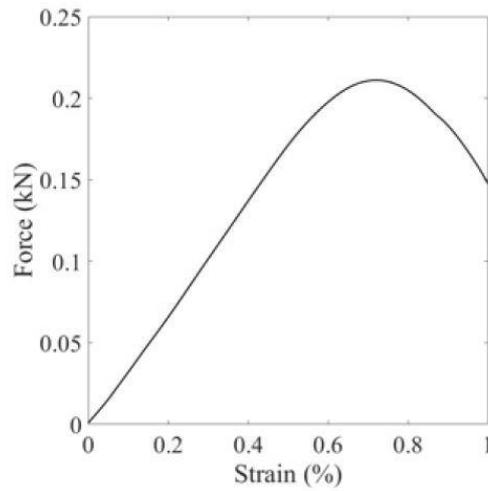

Figure 25: Loading curve on the top plate.

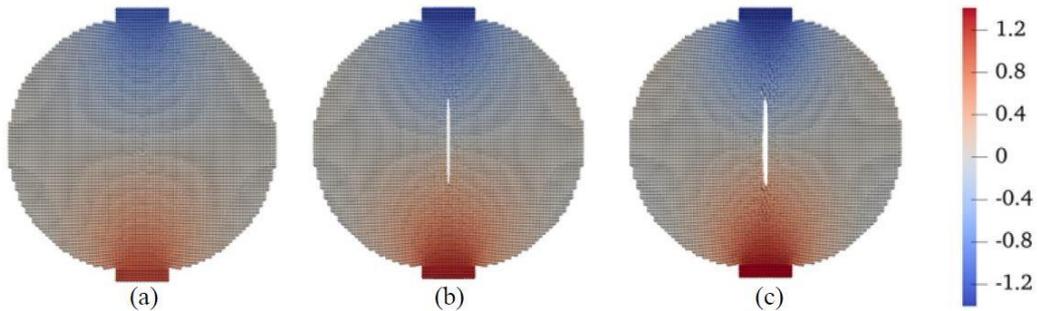

Figure 26: Contours of vertical displacement (mm) on deformed configuration at (a) $u = 1.1$ mm, (b) $u = 1.25$ mm, and (c) $u = 1.4$ mm.

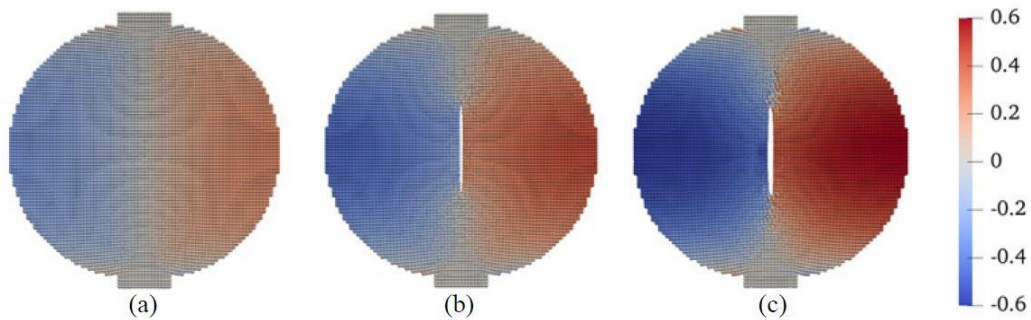

Figure 27: Contours of horizontal displacement (mm) on deformed configuration at (a) $u = 1.1$ mm, (b) $u = 1.25$ mm, and (c) $u = 1.4$ mm.

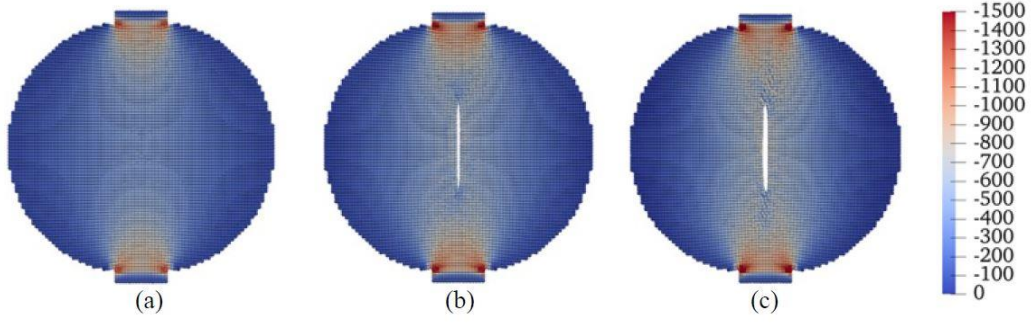

Figure 28: Contours of stress $\sigma_{yy}$ (kPa) on deformed configuration at (a) $u = 1.1$ mm, (b) $u = 1.25$ mm, and (c) $u = 1.4$ mm.

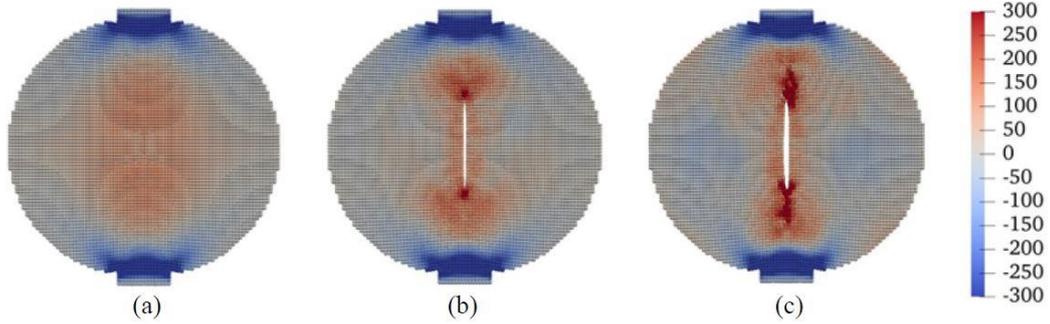

Figure 29: Contours of stress $\sigma_{xx}$ (kPa) on deformed configuration at (a) $u = 1.1$ mm, (b) $u = 1.25$ mm, and (c) $u = 1.4$ mm.

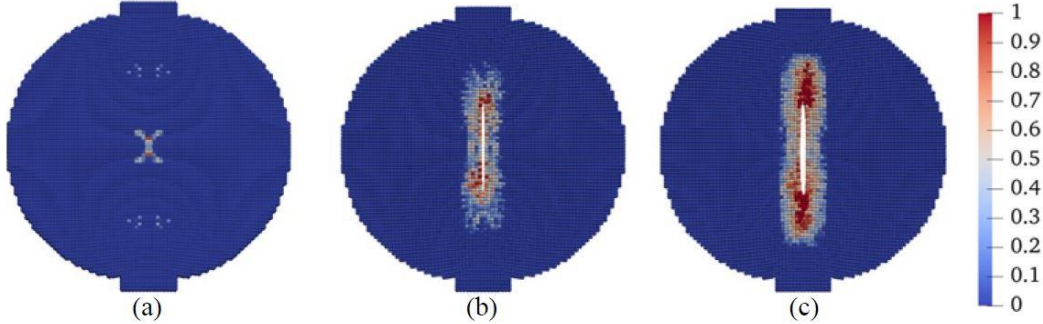

Figure 30: Contours of damage parameter on deformed configuration at (a) $u = 1.1$ mm, (b) $u = 1.25$ mm, and (c) $u = 1.4$ mm.

Second, in the base simulation we also investigate the influence of spatial discretization on the results by employing two distinct spatial discretization schemes. The first scheme utilizes 6224 points with a grid spacing (tlx) of 3.3 mm (referred to as grid 1), while the second scheme involves 10408 points with tlx = 2.5 mm (referred to as grid 2). Both schemes adopt the same horizon size, δ = 10 mm. The comparative results of these two discretization schemes are presented in Figures 31 through 36. Figure 31 displays the vertical loading curves obtained from simulations using both spatial discretization schemes. For a detailed comparison, we present contours of displacement and stress at a uniform displacement load (u = 1.25 mm) for both simulations as follows. Figure 32 compares the vertical displacement contours. Figure 33 contrasts the horizontal displacement contours. Figure 34 showcases the differences in vertical stress contours ($\sigma_{yy}$). Figure 35 illustrates the comparison in horizontal stress contours ($\sigma_{xx}$). Figure 36 compares the damage parameter contours at the same displacement load for both discretization schemes. The results from these comparisons suggest that, given a consistent horizon size, the choice of spatial discretization scheme exerts a minimal influence on the crack formation in the disk specimen under vertical compression loading.

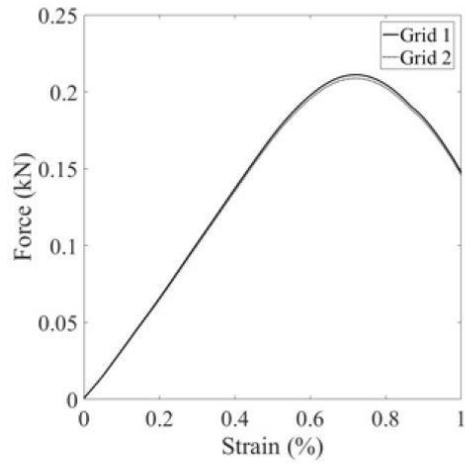

Figure 31: Comparison of the loading curve from the simulations with two spatial discretization schemes.

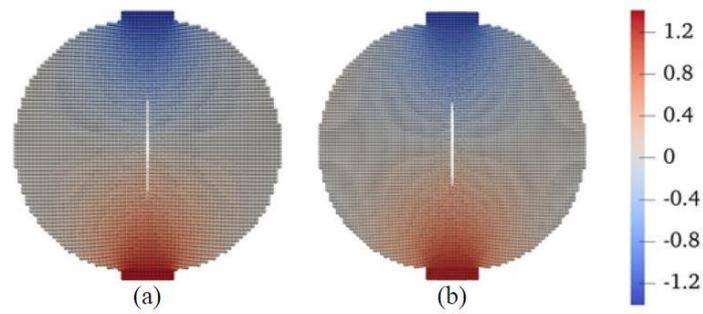

Figure 32: Contours of vertical displacement (mm) on deformed configuration at $u = 1.25$ mm: (a) grid 1 and (b) grid 2.

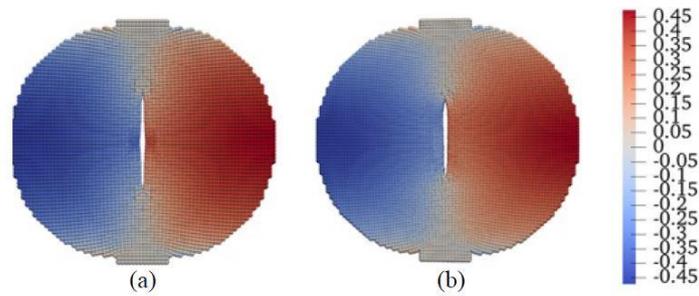

Figure 33: Contours of horizontal displacement (mm) on deformed configuration at $u = 1.25$ mm: (a) grid 1 and (b) grid 2.

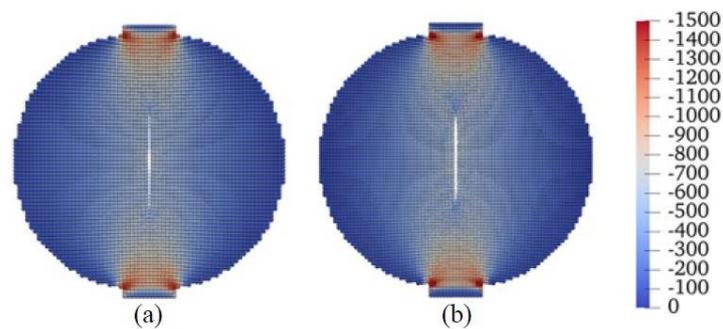

Figure 34: Contours of stress $\sigma_{yy}$ (kPa) on deformed configuration at $u = 1.25$ mm: (a) grid 1 and (b) grid 2.

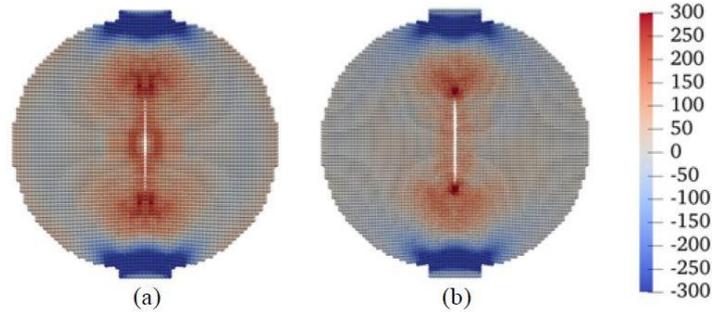

Figure 35: Contours of stress $\sigma_{xx}$ (kPa) on deformed configuration at $u = 1.25$ mm for (a) grid 1 and (b) grid 2.

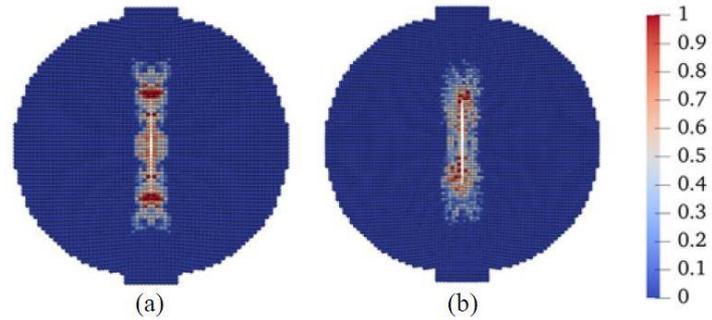

Figure 36: Contours of damage parameter on deformed configuration at $u = 1.25$ mm: (a) grid 1 and (b) grid 2.

Third, we investigate the effect of temperature variations on the cracking behavior in unsaturated elastic porous materials, specifically at 0ºC, 25ºC, and 50ºC. The other simulation parameters, including loading and spatial discretization, remain identical to those in the base simulation. The outcomes of this investigation are showcased in Figures 37 to 42. Figure 37 compares the loading curves from the simulations conducted at the three different temperatures, revealing a decrease in the peak load of the disk specimen with increasing temperature. The contours of vertical displacement at a uniform displacement load (u = 1.25 mm) across the three simulations are compared in Figure 38, while Figure 39 does the same for horizontal displacement. These figures illustrate that the rise in temperature increases the horizontal displacement, leading to a longer crack, whereas the vertical displacement is comparatively less influenced by temperature changes under the same loading conditions. The contours of vertical and horizontal stresses (σ-$_{yy}$ and σ-$_{xx}$, respectively) at u = 1.25 mm for the three simulations are presented in Figures 40 and

41. Finally, Figure 42 compares the contours of the damage variable at the same displacement load. Collectively, the results from Figures 40 to 42 indicate that the combination of temperature increase and mechanical loading significantly impacts the timing and progression of crack initiation and development in the specimen

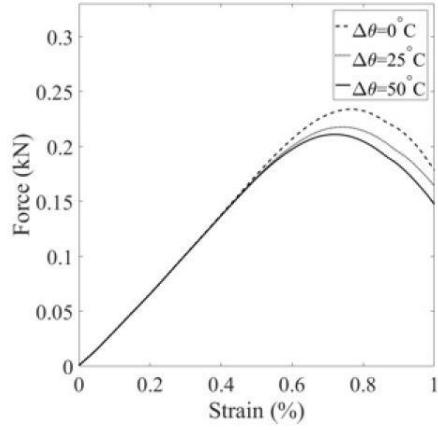

Figure 37: Comparison of the loading curves from the simulation with $\Delta\theta = 0°C$, $\Delta\theta = 25°C$, and $\Delta\theta = 50°C$.

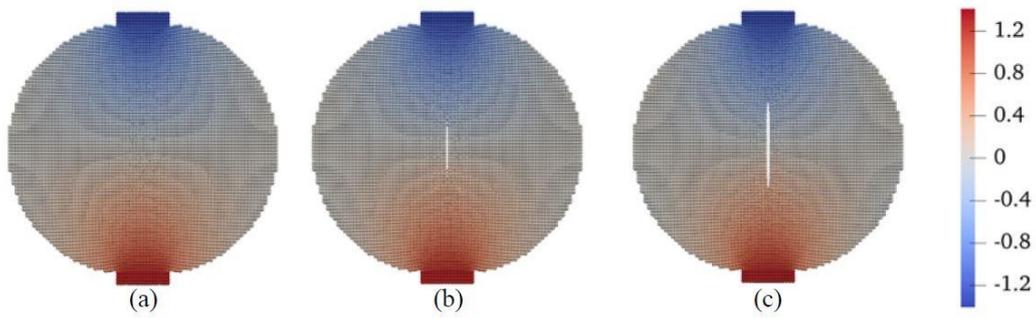

Figure 38: Contours of vertical displacement (mm) on deformed configuration at $u = 1.25$ mm: (a) $\Delta\theta = 0°C$, (b) $\Delta\theta = 25°C$, and (c) $\Delta\theta = 50°C$.

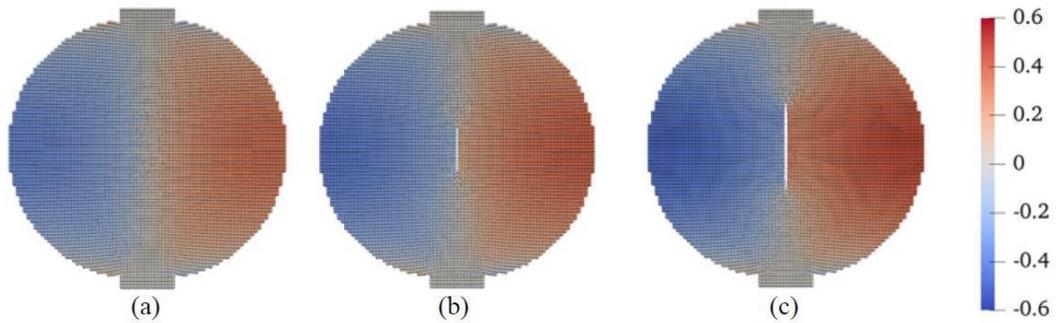

Figure 39: Contours of horizontal displacement (mm) on deformed configuration at $u = 1.25$ mm: (a) $\Delta\theta = 0°C$, (b) $\Delta\theta = 25°C$, and (c) $\Delta\theta = 50°C$.

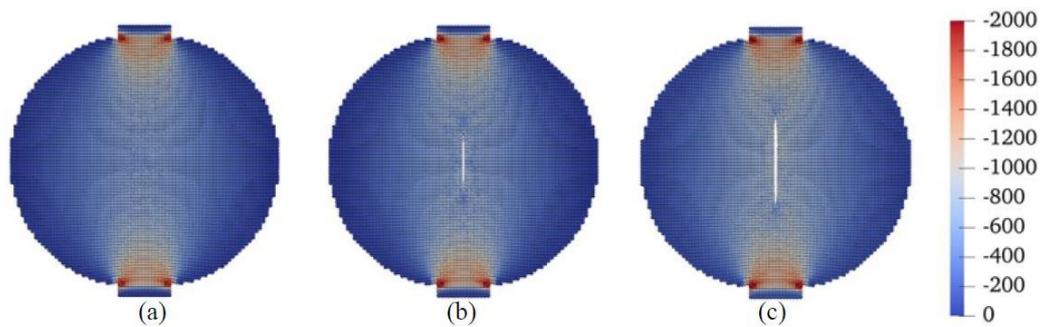

Figure 40: Contours of stress $\sigma_{yy}$ (kPa) on deformed configuration at $u = 1.25$ mm: (a) $\Delta\theta = 0°C$, (b) $\Delta\theta = 25°C$, and (c) $\Delta\theta = 50°C$.

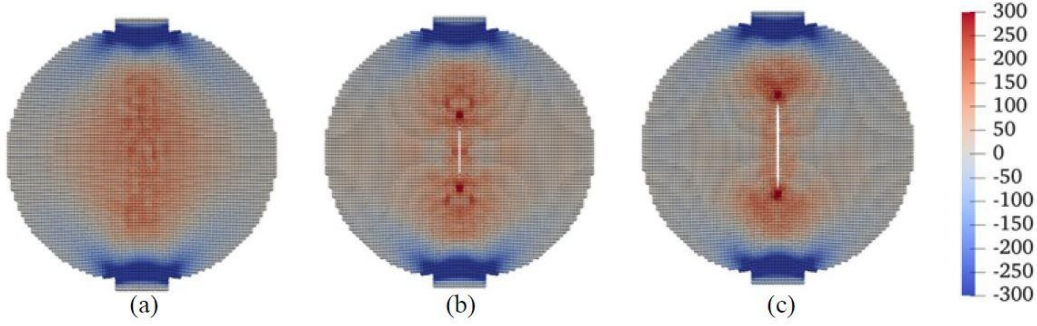

Figure 41: Contours of stress $\sigma_{xx}$ (kPa) on deformed configuration at $u = 1.25$ mm: (a) $\Delta\theta = 0°C$, (b) $\Delta\theta = 25°C$, and (c) $\Delta\theta = 50°C$.

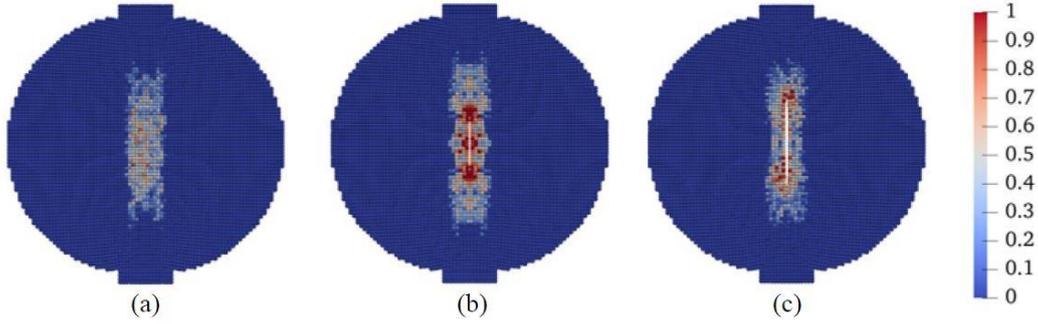

Figure 42: Contours of damage variable on deformed configuration at $u = 1.25$ mm: (a) $\Delta\theta = 0°C$, (b) $\Delta\theta = 25°C$, and (c) $\Delta\theta = 50°C$.

## 5. Summary

In this article, our investigation centers on shear banding and cracking in unsaturated porous media under non-isothermal conditions, utilizing a THM framework within PPM. A significant advancement of this study is the development of a nonlocal THM constitutive model specifically designed for unsaturated porous media in the PPM context. We have implemented the THM paradigm using an explicit Lagrangian meshfree algorithm, complemented by a return mapping algorithm for the numerical implementation of the nonlocal THM constitutive model. We have evaluated the performance and applicability of our proposed THM meshfree paradigm through a series of numerical examples. The results from these examples demonstrate the effectiveness and reliability of our THM PPM approach in accurately modeling the complex behaviors of shear banding and cracking in unsaturated porous media under varying thermal conditions. Importantly, our findings provide deep insights into the sophisticated relationship between temperature changes and the development of shear bands and cracks in such porous media under THM loading, highlighting the nuanced interdependencies in these phenomena.

## Acknowledgment

This work has been supported by the US National Science Foundation under contract number 1944009. The support is gratefully acknowledged. Any opinions or positions expressed in this article are those of the authors only and do not reflect any opinions or positions of the NSF.

## Conflict of Interest Statement

The authors declare no potential conflict of interest.

## Data Availability Statement

The data that support the findings of this study are available from the corresponding author upon reasonable request.